\documentclass[12pt]{amsart}

\usepackage{amsmath,amssymb,epic,lscape}

\theoremstyle{definition}

\theoremstyle{remark}

\numberwithin{equation}{section}

\def\DJ{{\hbox{D\kern-.8em\raise.15ex\hbox{--}\kern.35em}}}
\def\DJo{$\;$\kern-.4em
    \hbox{D\kern-.8em\raise.15ex\hbox{--}\kern.35em okovi\'c}}

\def\al{{\alpha}}

\def\la{{\lambda}}

\renewcommand{\subjclassname}{\textup{2000} Mathematics Subject
Classification }

\begin{document}

\title[Hadamard matrices of order $764$]
{Hadamard matrices of order $764$ exist}

\author[D.\v{Z}. \DJ okovi\'{c}]
{Dragomir \v{Z}. \DJ okovi\'{c}}

\address{Department of Pure Mathematics, University of Waterloo,
Waterloo, Ontario, N2L 3G1, Canada}

\email{djokovic@uwaterloo.ca}

\thanks{
The author was supported by the NSERC Grant A-5285.}


\date{}

\begin{abstract}
Two Hadamard matrices of order $764$ of Goethals--Seidel type
are constructed.
\end{abstract}

\maketitle
\subjclassname{ 05B20, 05B30 }
\vskip5mm

Recall that a Hadamard matrix of order $m$ is a $\{\pm1\}$-matrix
$A$ of size $m\times m$ such that $AA^T=mI_m$, where
$T$ denotes the transpose and $I_m$ the identity matrix.
We refer the reader to one of \cite{KH,SY} for the survey of
known results about Hadamard matrices.

In our previous note \cite{DZ}, written about 13 years ago,
we listed 17 integers $n\le500$ for which no Hadamard matrix
of order $4n$ was known at that time. Two of these integers were
removed in that note and the smallest one, $n=107$, was
removed recently by Kharaghani and Tayfeh-Rezaie \cite{KT}.
Among the remaining 14 integers $n$ only four are less than
1000. The problem of existence of Hadamard matrices of
these four orders, namely 668, 716, 764 and 892, has been
singled out as Research Problem 7 in the recent book \cite{KH}
by Kathy Horadam. In this note we shall remove the integer 764
from the mentioned list by constructing two examples of Hadamard
matrices of Goethals--Seidel type of that order.
(We have constructed a bunch of examples but we shall present
only two of them.)
Consequently, the revised list now consists of 13 integers:
\[ 167,\, 179,\, 223,\, 251,\, 283,\, 311,\, 347,\, 359,\,
419,\, 443,\, 479,\, 487,\, 491; \]
all of them primes congruent to $3 \pmod{4}$.

For the remainder of this note we set $n=191$. 
Let $G$ be the multiplicative group of non-zero residue classes
modulo the prime $n=191$, a cyclic group of order $n-1=190$, and let
$H=\langle 39 \rangle =\{ 1,39,184,109,49 \}$ be its subgroup
of order 5. We choose the enumeration of the 38 cosets $\al_i$,
$0\le i\le 37$, of $H$ in $G$ so that $\al_{2i+1}=-1\cdot\al_{2i}$
for $0\le i\le 18$ and
\[
\begin{array}{lllll}
\al_0=H, \quad & \al_2=2H, \quad & \al_4=3H, \quad & \al_6=4H, \quad & \al_8=6H, \\
\al_{10}=8H, & \al_{12}=9H, & \al_{14}=11H, & \al_{16}=12H, & \al_{18}=13H, \\
\al_{20}=16H, & \al_{22}=17H, & \al_{24}=18H, & \al_{26}=19H, & \al_{28}=22H, \\
\al_{30}=32H, & \al_{32}=36H, & \al_{34}=38H, & \al_{36}=41H. &
\end{array}
\]

Define four index sets:
\begin{eqnarray*}
J_1 &=& \{ 1,7,9,10,11,13,17,18,25,26,30,31,33,34,35,36,37 \}, \\
J_2 &=& \{ 1,4,7,9,11,12,13,14,19,21,22,23,24,25,26,29,36,37 \}, \\
J_3 &=& \{ 0,3,4,5,7,8,9,16,17,19,24,25,29,30,31,33,35,37 \}, \\
J_4 &=& \{ 1,3,4,5,8,11,14,18,19,20,21,23,24,25,28,29,30,32,34,35 \}
\end{eqnarray*}
and introduce the following four sets of integers modulo 191:
\[ S_k = \bigcup_{i\in J_k} \al_i,\quad k=1,2,3,4. \]
Their cardinals $n_k=|S_k|=5|J_k|$ are:
\[ n_1=85,\, n_2=n_3=90,\, n_4=100 \]
and we set
\[ \la = n_1+n_2+n_3+n_4-n=174. \]

For $r\in\{ 1,2,\ldots,190 \}$ let $\la_k(r)$ denote the number of
solutions of the congruence $i-j\equiv r \pmod{191}$ with
$\{i,j\}\subseteq S_k$. It is easy to verify (by using a computer) that
\[ \la_1(r)+\la_2(r)+\la_3(r)+\la_4(r)=\la \]
is valid for all such $r$. Hence the sets
$S_1,S_2,S_3,S_4$ are supplementary
difference sets (SDS), with associated decomposition
\begin{eqnarray*}
4n=764 &=& 9^2+11^2+11^2+21^2 \\
&=& \sum_{k=1}^4 (n-2n_k)^2.
\end{eqnarray*}

Let $A_k$ be the $n\times n$ circulant matrix with first row
\[ a_{k,0},\, a_{k,1},\, \ldots ,\, a_{k,n-1} \]
where $a_{k,j}=-1$ if $j\in S_k$ and $a_{k,j}=1$ otherwise.
These $\{\pm1\}$-matrices satisfy the identity
\[ \sum_{k=1}^4 A_kA_k^T=4nI_n. \]

One can now plug in the matrices $A_k$ into the Goethals--Seidel
array to obtain a Hadamard matrix of order $4n=764$.

Our second example is constructed in the same way by using
the index sets:
\begin{eqnarray*}
J_1 &=& \{ 0,1,6,8,9,11,12,16,18,20,21,23,28,31,33,36,37 \}, \\
J_2 &=& \{ 0,1,3,4,10,12,13,17,20,22,24,31,32,33,34,35,36,37 \}, \\
J_3 &=& \{ 4,8,9,10,12,13,14,16,17,20,21,24,26,27,29,31,32,34 \}, \\
J_4 &=& \{ 1,7,9,10,11,12,14,15,16,17,20,22,23,25,28,29,32,33,34, \\
&& \quad 37 \}.
\end{eqnarray*}
The two solutions are not equivalent in the sense that the
two SDS's are not equivalent. (For the definition of equivalence
for SDS's see our note \cite{DZ}.)

\end{document}